\documentclass{siamltex} 

\usepackage{mathtools}  % amsmath + more
\usepackage{amssymb}    % More symbols
\usepackage{amstext}    % Text inside mathenv
\usepackage{commath}
\usepackage{amsmath}

\usepackage[utf8]{inputenc}
\usepackage[T1]{fontenc}
\usepackage{microtype}

\usepackage{algorithm}
\usepackage{algpseudocode}

\usepackage{fvextra}        % Req by minted (must load before csquotes)
\usepackage{csquotes}       % For bibliography quotations
\usepackage[dvipsnames]{xcolor}
\usepackage{color, xcolor}     % more colors
\usepackage{tcolorbox}

\usepackage{soul}

\usepackage{amsfonts}
\usepackage{amssymb}
\usepackage{amsmath}
\usepackage{graphicx} 
\graphicspath{ {./figures/} }

\usepackage{comment}
\usepackage{varioref} % reference far away (1)
\usepackage[colorlinks = true]{hyperref} % links in references (2)
\usepackage{cleveref} % smart references (3)
\hypersetup{
   linkcolor={black},
   citecolor={black},
   urlcolor={blue!80!black}
}
\usepackage{url}

\tcbset{width=\textwidth,boxrule=0pt,colback=yellow,arc=0pt,auto outer arc,left=0pt,right=0pt,boxsep=5pt}

\usepackage[bottom]{footmisc} % Footnotes at bottom of page

\usepackage{float}               % Force tables and images position (H and H!)
\usepackage{wrapfig}            % Wrap images like in HTML

\usepackage{caption}
\usepackage{subcaption}                     % Subfigures
\usepackage[rightcaption]{sidecap}
\usepackage{xltabular, colortbl, booktabs} % longtable + tabularx (has bug with booktabs: fix below)

\usepackage{array, multirow, multicol, makecell}

\makeatletter
\def\@BTrule[#1]{%
  \ifx\longtable\undefined
    \let\@BTswitch\@BTnormal
  \else\ifx\hline\LT@hline
    \nobreak
    \let\@BTswitch\@BLTrule
  \else
     \let\@BTswitch\@BTnormal
  \fi\fi
  \global\@thisrulewidth=#1\relax
  \ifnum\@thisruleclass=\tw@\vskip\@aboverulesep\else
  \ifnum\@lastruleclass=\z@\vskip\@aboverulesep\else
  \ifnum\@lastruleclass=\@ne\vskip\doublerulesep\fi\fi\fi
  \@BTswitch}
\makeatother
\let\norm\undefined % <-- "Undefine" \norm
\DeclarePairedDelimiter\norm{\lVert}{\rVert}
\def\cvd{~\vbox{\hrule\hbox{%
   \vrule height1.3ex\hskip0.8ex\vrule}\hrule } }

\author{Andrew J. Higgins\thanks{Los Alamos National Laboratories, Los Alamos, New Mexico, USA (\texttt{ahiggins@lanl.gov}).} \and Logan O. Reed\thanks{Temple University, Philadelphia, Pennsylvania, USA (\texttt{lreed@temple.edu}, \texttt{szyld@temple.edu}).} \and Daniel B. Szyld\footnotemark[3]}

\date{}

\title{Near-Optimal Convergence of Widlund's method and Numerical Exploration of Different Implementations%
    \thanks{this version dated October 29, 2026}.}

\begin{document}
%%%%%%%%%%%%%%%%%%%%%%%%%%%%%%%%%%%%%%%%%%%%%%%%%%%%%%%%%%%%%%%%%%%%%%%%%%%%%%%%
% TITLE
%%%%%%%%%%%%%%%%%%%%%%%%%%%%%%%%%%%%%%%%%%%%%%%%%%%%%%%%%%%%%%%%%%%%%%%%%%%%%%%%

    % Default title or use titlepage.tex

      % the following is the abstract in one line for application to conferences :)
      % For linear systems arising from the discretization of certain classes of mathematical models, the splitting of the linear system into its Hermitian and skew-Hermitian parts has a physical meaning. If, in addition, the Hermitian part is positive definite, there exist Krylov subspace methods based on three-term recurrences where the Hermitian part is used as a preconditioner. Two such Krylov subspace methods are Widlund's method and Rapoport's method. In this paper we establish a ``near-optimal'' guarantee for the convergence of Widlund's method, showing that the overall convergence of Widlund's method is nearly as good as that of Rapoport's method. Specifically, we show that for every iteration $k$, there exists an iteration $j\leq k$ where the residual norm of Widlund's method is less than $\sqrt{k+1} $ times the residual norm of Rapoport's method at iteration $k$. This bound is sharp. We also provide a computational study on various implementations of Widlund's method in the exact and inexact cases, and show that certain implementations may stagnate while others do not, i.e., we look at the robustness of the different implementations. We recommend two implementations.

    \maketitle

    \begin{abstract}
      For linear systems arising from the discretization of certain classes of mathematical models, the splitting of the linear system into its Hermitian and skew-Hermitian parts has a physical meaning. 
      If, in addition, the Hermitian part is positive definite, there exist Krylov subspace methods based on three-term recurrences where the Hermitian part is used as a preconditioner. 
      Two such Krylov subspace methods are Widlund's method and Rapoport's method.
      In this paper we establish a ``near-optimal'' guarantee for the convergence of Widlund's method, showing that the overall convergence of Widlund's method is nearly as good as that of Rapoport's method.
      Specifically, we show that for every iteration $k$, there exists an iteration $j\leq k$ where the residual norm of Widlund's method is less than $\sqrt{k+1} $ times the residual norm of Rapoport's method at iteration $k$, and this bound is sharp.
      We also provide a computational study on various implementations of Widlund's method in the exact and inexact cases, and show that certain implementations may stagnate while others do not, i.e., we look at the robustness of the different implementations. 
      We recommend two implementations.
    \end{abstract}


%%%%%%%%%%%%%%%%%%%%%%%%%%%%%%%%%%%%%%%%%%%%%%%%%%%%%%%%%%%%%%%%%%%%%%%%%%%%%%%%
% TOC & lists
%%%%%%%%%%%%%%%%%%%%%%%%%%%%%%%%%%%%%%%%%%%%%%%%%%%%%%%%%%%%%%%%%%%%%%%%%%%%%%%%

    \pagenumbering{Roman}

    %\setcounter{tocdepth}{2}
    % \tableofcontents \pagebreak

    % \listoffigures \pagebreak
    % \listoftables \clearpage

    \pagenumbering{arabic}

%%%%%%%%%%%%%%%%%%%%%%%%%%%%%%%%%%%%%%%%%%%%%%%%%%%%%%%%%%%%%%%%%%%%%%%%%%%%%%%%
% SECTIONS
%%%%%%%%%%%%%%%%%%%%%%%%%%%%%%%%%%%%%%%%%%%%%%%%%%%%%%%%%%%%%%%%%%%%%%%%%%%%%%%%

    % Paragraph spacing (placed after ToC)
    \setlength{\parskip}{1em plus 0.5em minus 0.2em}
    %\setlength{\parindent}{0pt}

    %! TEX root = **/000-main.tex
% vim: spell spelllang=en:

\section{Introduction}%
\label{sec:introduction}

We want to solve the equation
\begin{equation}
  Ax = b \label{LinearEquation}
,\end{equation}
where $A \in \mathbb{C}^{n \times n}$. 
One can split $A$ into its Hermitian and skew-Hermitian components, $H$ and $S$, i.e., $A = H + S$, where
\begin{equation}
    H = \frac{1}{2}\left( A + A^* \right), \quad S = \frac{1}{2}\left( A - A^* \right) \label{HSformula}
.\end{equation}

We assume $H$ to be positive definite, and so we can define the $H$-inner product $\langle x,y \rangle_H = y^*Hx$, which induces a norm $\norm[\big]{x}_H^2 = x^*Hx = \norm{H^{1/2}x}^2$. 
For reference, note that $\norm{Hx}^2 = (Hx)^*Hx = x^*H^2x = \norm{x}_{H^2}^2$ and $\|H^{-1}x\|_{H} = \|x\|_{H^{-1}}$.

For linear systems arising from the discretization of certain classes of mathematical models, the splitting $A = H + S$ arises naturally and has a physical meaning. 
This occurs, e.g., in the discretization of dissipative Hamiltonian DAEs and of port-Hamiltonian systems \cite{flexible-short-recurrence-Frommer, Guducu.Liesen.Mehrmann.Szyld.22}.
In these settings $H$ is positive definite, and there exist Krylov subspace methods based on three-term recurrences where the Hermitian part is used as a preconditioner.

Define $K = H^{-1}S$.
Let $I_n$ be the $n \times n$ identity matrix, or simply $I$ when the order is understood from the context.
We can re-formulate $A$ as $A = H(I+K)$.
It follows that (\ref{LinearEquation}) is equivalent to
\begin{equation}
    (I_n + K)x = \hat{b} \label{IKsystem}
,\end{equation}
where $\hat{b} = H^{-1}b$.

Because $I + K$ is a shifted skew-Hermitian matrix with respect to the $H$-inner product, it has an optimal 3-term recurrence for the construction of an orthogonal basis of a corresponding Krylov subspace, see, e.g., \cite{recurrence-faber-manteuffel,short-recurrence-liesen-strakos}.
Using the notation $\mathcal{K}_{m}\left( A, b \right)$ for the Krylov subspace $\mathcal{K}_{m}\left( A, b \right) = \text{span} \left\{ b,Ab,\ldots,A^{m-1}b \right\}$ of dimension $m$ generated by the matrix $A$ and vector $b$, the 3-term recurrence can be written as
 \begin{equation}
  \left( I + K \right)V_{m} = V_{m+1}T_{m+1,m} \label{eq:lanczosRelation}
,\end{equation}
where the columns $v_{i}$ of the $n \times m$ matrix $V_{m}$ form an $H$-orthonormal basis of \linebreak$\mathcal{K}_{m}( I+K, \hat{b} )$, and the $(m+1)\times m$ matrix $T_{m+1,m}$ is tridiagonal.
Note that \linebreak $\mathcal{K}_m(I + K, \hat{b}) = \mathcal{K}_m (K,\hat{b})$, as Krylov subspaces are shift invariant.
Let $T_{m}$ be the matrix consisting of the first $m$ rows of  $T_{m+1,m}$.

The method of Widlund \cite{Widlund78} and, independently, of Concus and Golub \cite{ConcusGolub} is an oblique projection method onto this Krylov subspace, $\mathcal{K}_{m}( K, \hat{b} )$, designed so that
\begin{equation}
  x_m \in x_0 + \mathcal{K}_m(K, \hat{b}), \quad r_m = \hat{b} - \left( I + K \right)x_m \perp_{H} \mathcal{K}_m(K, \hat{b})
    \label{eq:CGWkrylovCondtion}
,\end{equation}
where $x_0$ is the initial vector.
This is an example of an Orthogonal Residual (OR) method onto $\mathcal{K}_{m}( K, \hat{b} )$, as defined, e.g., in \cite{Eiermann_Ernst_2001}.

The method of Rapoport \cite{nonlinear-lanzcos-rapoport} is the corresponding Minimum Residual (MR) method using the same Krylov subspace, $\mathcal{K}_{m}( K, \hat{b} )$.
That is, Rapoport's method is an orthogonal projection method designed so that 
\begin{equation}
  \left\| r_{m} \right\|_{H} = \min_{x \in x_{0} + \mathcal{K}_{m}( K, \hat{b} )} \| \hat{b} -  \left( I + K \right)x \|_{H} \label{eq:rapoportVariationalDefinition}
.\end{equation}
Throughout the paper, we will use $r_{k}^{W}$ and $r_{k}^{R}$ to denote the $k$-th residual of Widlund's method and Rapoport's method, respectively.
We remark that due to the methods of Widlund/Rapoport being an OR/MR pair on $\mathcal{K}_{m}( K, \hat{b} )$, the pair can be seen as FOM/MINRES on the preconditioned system \eqref{IKsystem} with the $H$-inner product replacing the standard inner product.

In \cite{Guducu.Liesen.Mehrmann.Szyld.22}, both Widlund's and Rapoport's methods are considered, in part because Rapoport's method, as an MR method, has monotonically decreasing residual norms. 
Our first contribution in this paper is to show that the residual norms of Widlund's method are ``nearly-optimal'' in a certain sense.
Thus, the benefit of Rapoport's monotonically decreasing residual norms may not outweigh the associated computational cost, see, e.g., the timings reported in \cite{Guducu.Liesen.Mehrmann.Szyld.22}.
For this reason, we will emphasize Widlund's method from Section \ref{sec:currentsituation} onward.

Our second contribution in this paper is to provide a computational study demonstrating that different implementations of Widlund's method, while mathematically equivalent, may have different behaviors.
We consider three implementations based on~\cite{Widlund78}, one straightforward implementation of \eqref{eq:lanczosRelation} used in \cite{MehrmannStoll}, and one which uses right preconditioning to create a flexible method from \cite{flexible-short-recurrence-Frommer}.
We also consider the practical situation in which the matrix-vector product with $K = H^{-1}S$ is not applied exactly.
In other words, in order to compute $Kv$ we first compute $Sv$, then solve a linear system $Hw = Sv$. 
If this system is solved iteratively or inexactly, then the basis vectors of $\mathcal{K}_{m}$ may fail to be $H$-orthogonal.
Thus, the resulting ``inexact Widlund's method'' is no longer guaranteed to have an orthogonal basis, and therefore the method is also non-optimal in the sense of \cite{non-opt}. 
Again, there are two situations here: the inexact matrix-vector products and, as a consequence of this inexactness, the matrices $V_m$ which are no longer $H$-orthogonal.
We assume exact arithmetic throughout; that is, the word ``inexact" refers only to the matrix-vector product, and does not refer to round-off error.
We show that not all implementations appear to be stable in the inexact setting.

The paper is organized as follows.
In Section \ref{sec:nearOptimal} we show that the residual norms of Widlund's method are ``nearly-optimal'' in the sense of \cite{near-optimal-chen-meurant}.
In Section \ref{sec:currentsituation} we provide a computational study of different implementations of Widlund's method, and in Section~\ref{sec:inexAlgs} we provide a computational study of the implementations when $H^{-1}$ is computed inexactly, thus discussing the robustness of the different implementations.

    %! TEX root = **/000-main.tex
% vim: spell spelllang=en:

\section{Near Optimal Convergence of Widlund}%
\label{sec:nearOptimal}

In this section we aim to show that Widlund's method is  ``near-optimal'' in the sense of \cite{near-optimal-chen-meurant}.
We take the same initial vector, $x_{0}^{W} = x_{0}^{R}$, for both methods, and thus $r_{0}^{W} = r_{0}^{R}$.
First of all, since Rapoport's method is an MR method, $\| r_{k}^{R} \|_{H} \leq \| r_{k}^{W} \|_{H}$ for all $k$, and Widlund's method will never converge in fewer iterations than Rapoport's method.
However, Widlund's method is ``near-optimal'' in the sense that, on average, $\|  r_{k}^{W}\|_{H}$ is no worse than $\sqrt{k+1} \| r_{k}^{R} \|_{H}$.

\begin{theorem} \label{thm:nearOptimal}
  For every $k \geq 1$,
  \begin{equation}
    \min_{0\leq j \leq k} \| r_{j}^{W} \|_{H} \leq \sqrt{k+1} \| r_{k}^{R} \|_{H}
    \label{eq:nearOptimal}
  \cdot\end{equation}
\end{theorem}

Our proof mirrors the argument found in \cite[Theorem 2.1]{near-optimal-chen-meurant} regarding the near-optimality of FOM.
Central to this proof is a characterization of the Rapoport residuals in terms of the Widlund residuals.
Define $\vartheta_{1} = 1$, and for $k \geq 1$ define
\[
\vartheta_{k+1} = - \frac{1}{t_{k+1,k}}\sum_{j=k-1}^{k} \vartheta_{j} t_{j,k}
,\] 
where $t_{i,j}$ is the $\left( i,j \right)$-th entry of $T_{k+1,k}$.
The residuals are related in the following sense.

\begin{theorem} \label{thm:residualFromInv}
  For every $k \geq 0$,
  \begin{equation}
    \frac{\| r_{k}^{W} \|_{H}}{\| r_{0}^{W} \|_{H}}=\frac{1}{\left| \vartheta_{k+1} \right|}
    \label{eq:residualFromInvWidlund}
  ,\end{equation}
  \begin{equation}
    \frac{\| r_{k}^{R} \|_{H}^{2}}{\| r_{0}^{R} \|_{H}^{2}} = \left( \sum_{j=1}^{k+1} \left| \vartheta_{j} \right|^{2} \right)^{-1}
    \label{eq:residualFromInvRapoport}
  \cdot\end{equation}
\end{theorem}

This is a special case of \cite[Theorem 3.12]{nonsymmetric-krylov-meurant-tebbens}, which requires that the two methods are the orthogonal residual and minimum residual pair over a given Krylov subspace.
Because Widlund's and Rapoport's methods are the OR/MR pair on the space $\mathcal{K}(K, \hat{b})$, Theorem~\ref{thm:residualFromInv} holds.

{\em Proof of Theorem \ref{thm:nearOptimal}.}
By \cite[Lemma 2.3]{nonsymmetric-krylov-meurant-tebbens}, $T_{k}$ is singular if and only if \mbox{$\vartheta_{k+1}=0$}.
In this case, the residual of Widlund's method is undefined and we define \mbox{$ \|   r_{k+1}^{W} \|_{H} = \infty$}.
By Theorem \ref{thm:residualFromInv}, for all $0 \leq j \leq k$
\[
\frac{\| r_{0}^{W} \|_{H}}{\| r_{j}^{W} \|_{H}} = \left| \vartheta_{j+1} \right|
\cdot\] 

\noindent Plugging \eqref{eq:residualFromInvWidlund} into \eqref{eq:residualFromInvRapoport} we obtain
\[
  \frac{\| r_{k}^{R} \|_{H}^{2}}{\| r_{0}^{R} \|_{H}^{2}} = \left( \sum_{j=0}^{k} \frac{\| r_{0}^{W} \|_{H}^{2}}{\| r_{j}^{W} \|_{H}^{2}} \right)^{-1}
\cdot\] 

\noindent From our construction, $\| r_{0}^{R} \|_{H} = \| r_{0}^{W} \|_{H}$, and thus
\begin{equation}
\| r_{k}^{R} \|_{H}^{2} = \left( \sum_{j=0}^{k} \frac{1}{\| r_{j}^{W} \|_{H}^{2}} \right)^{-1}
\label{eq:intermediateNearOptimal}
\cdot\end{equation}

\noindent Taking the largest $\frac{1}{\| r_{j}^{W} \|_{H}^{2}}$,
\begin{align*}
  \| r_{k}^{R} \|_{H}^{2} &\geq \left( \left( k+1 \right) \max_{0\leq j\leq k} \frac{1}{\| r_{j}^{W} \|_{H}^{2}} \right)^{-1}\\
                             &= \left( k+1 \right)^{-1}\min_{0 \leq j \leq k}\| r_{j}^{W} \|_{H}^{2}
\cdot\end{align*}

\noindent Taking the square root of both sides, we obtain the desired inequality.
\cvd

We remark that we do not use any specific property of Widlund/Rapoport, so the proof of Theorem \ref{thm:nearOptimal} gives the same bound on any OR/MR pair.

\begin{figure}[htb]
    \centering
    \includegraphics[width=\textwidth]{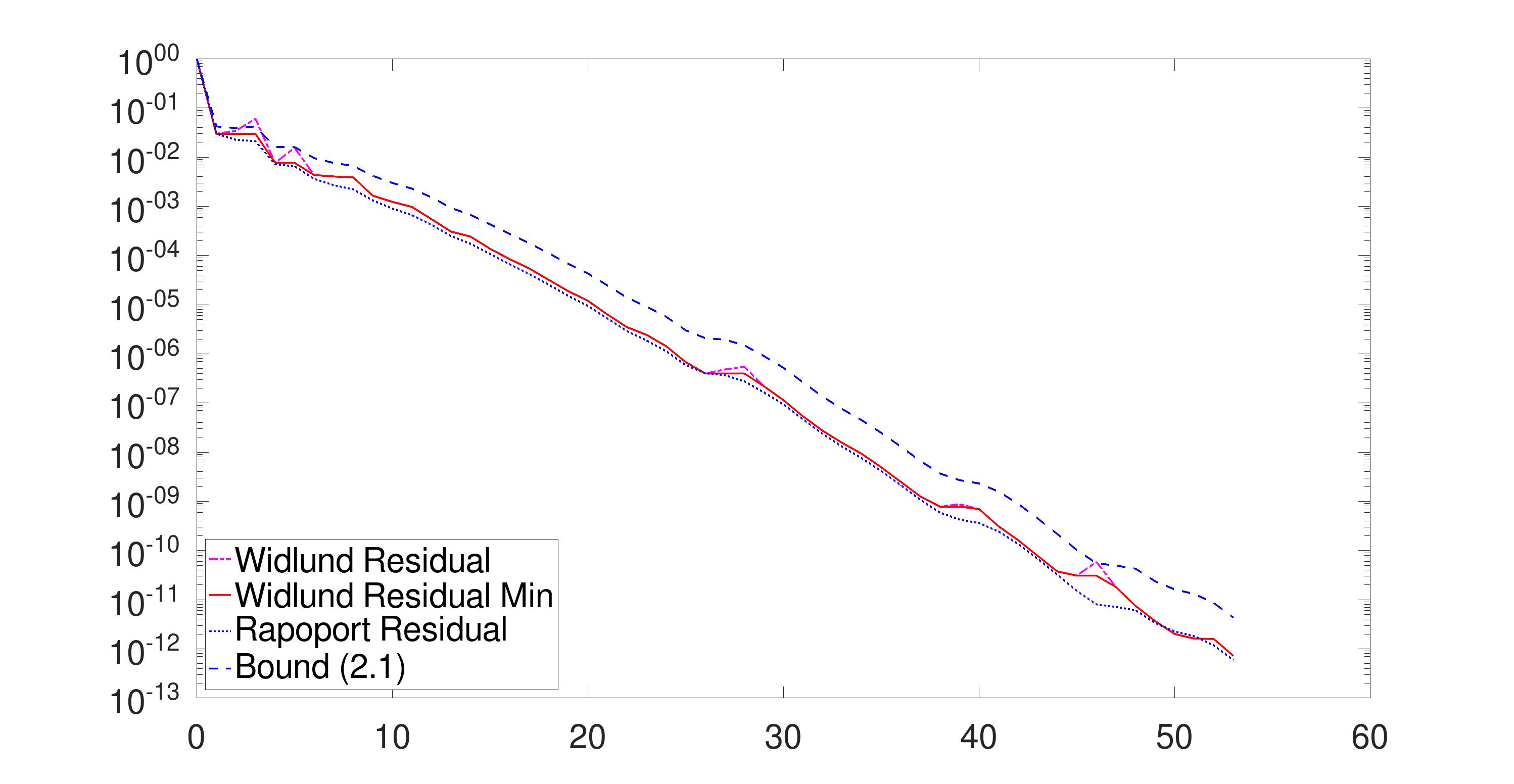}
    \caption{Left-hand and right-hand sides of \eqref{eq:nearOptimal} alongside the $H$-norm relative residual.}
    \label{fig:prob3NearOptimal}
\end{figure}

In Figure \ref{fig:prob3NearOptimal} we illustrate \eqref{eq:nearOptimal} holding for \cite[Example~2]{flexible-short-recurrence-Frommer} of order $n=3,151,875$, which was originally considered in \cite[Section 6.2]{Guducu.Liesen.Mehrmann.Szyld.22}.
It is the system created by discretizing the incompressible Stokes equations using the $Q_{1}$-$Q_{1}$ finite element approximation of the (unsteady) channel domain problem in IFISS \cite{ifiss}.

The bound \eqref{eq:nearOptimal} is tight, as we show in the following theorem.

\begin{theorem}\label{thm:nearOptimalTightness}
  For every $k \geq 1$, there exists a matrix $A$ and vector $b$ such that 
   \begin{equation}
    \min_{0\leq j \leq k} \| r_{j}^{W} \|_{H} = \sqrt{k+1} \| r_{k}^{R} \|_{H}
  \cdot\end{equation} 
\end{theorem}

{\em Proof.}
Fix $k\geq 1$. Let $e_j,$ $j=1,\dots,k+1$ be the standard unit vectors. Let $H=I_{k+1}$, $S=\mathrm{i}\Sigma$, where $\Sigma\in\mathbb{R}^{(k+1)\times(k+1)}$ is the symmetric tridiagonal matrix
\[
    \Sigma=
    \begin{bmatrix}
        0 & 1 \\
        1 & 0 & 2 \\
        & \ddots & \ddots & \ddots \\
        && k-1 & 0 & k \\
        &&& k & 0
    \end{bmatrix}.
\]

Let $A = H + S$ and $b = e_1$.
Since $H$ is symmetric positive definite and $S$ is skew-symmetric, we can define Widlund's method (with initial vector $x_0^W=0$) on $\mathcal{K}\left( I + H^{-1}S, H^{-1}b \right)$.
Let $j \leq k$.
Expanding the Arnoldi relation, we can obtain
\[
    T_j=
    \begin{bmatrix}
        1 & -1 \\
        1 & 1 & -2 \\
        & \ddots & \ddots & \ddots \\
        && j-2 & 1 & -(j-1) \\
        &&& j-1 & 1
    \end{bmatrix}.
\]
Since Widlund's method is an OR method, the residual is $\|r_j^W\|_H = j |e_j^T y_j|$, where $y_j$ is defined such that $T_j y_j = e_1$. 
We seek an explicit expression for $\|r_j^W\|_H$, which we obtain using the cofactor formula for $T_j^{-1}$.
Let $d_1 = 1$, $d_2 = 2$, and $d_i = \text{det}(T_i)$ for all $3\leq i\leq j$.
Since $T_j$ is tridiagonal, we can determine $d_j$ using the recurrence relation
\[
  d_j = d_{j-1} + (j-1)^2 d_{j-2}.
\]
From this, one can show that $d_j = j!$.
Hence $T_j$ is nonsingular, and the cofactor formula gives 
\[
  e_j^T y_j = (T_j^{-1})_{j,1} = (-1)^{j+1}\frac{\prod_{\ell = 1}^{j-1} \ell}{d_j}.
\]
Thus, the residual is 
\[
  \| r_j^W \|_H = j |e_j^T y_j| = \frac{\prod_{\ell = 1}^{j} \ell}{d_j} = \frac{j!}{j!} = 1.
\]
Further, since $H = I_{k+1}$, $\|r_0^W\|_H = \|e_1\|_2 = 1$ and so $\|r_j^W\|_H = 1$ for all $j = 0,\ldots,k$.
From \eqref{eq:intermediateNearOptimal}, we have
\[
  \|r_k^R\|_H^2 = \left( \sum_{j=0}^{k} \frac{1}{\| r_{j}^{W} \|_{H}^{2}} \right)^{-1} = (k+1)^{-1}.
\]
Finally, since $\min_{0\leq j \leq k}\|r_j^W\|_H = 1$, $\min_{0\leq j \leq k}\|r_j^W\|_H = \sqrt{k+1} \|r_k^R\|_H$.\cvd

We remark that \eqref{eq:intermediateNearOptimal} can be used to show the following well-known relationship (see, e.g., \cite[\S4]{krylov-survey-simoncini-szyld})
\[
\| r_{k}^{W} \|_{H} = \frac{\| r_{k}^{R} \|_{H}}{\sqrt{1 - \left(\frac{\| r_{k}^{R} \|_{H}}{\| r_{k-1}^{R} \|_{H}}\right)^{2}} }
\cdot\] 
Indeed, from \eqref{eq:intermediateNearOptimal},
\[
\frac{1}{\| r_{k}^{W} \|_{H}^{2}} = \sum_{j=0}^{k} \frac{1}{\| r_{j}^{W} \|_{H}^{2}} - \sum_{j=0}^{k-1} \frac{1}{\| r_{j}^{W} \|_{H}^{2}} = \frac{1}{\| r_{k}^{R} \|_{H}^{2}} - \frac{1}{\| r_{k-1}^{R} \|_{H}^{2}}
\cdot\] 
Then, taking the reciprocal of both sides,
\begin{align*}
  \| r_{k}^{W} \|_{H}^{2} &= \frac{1}{\frac{1}{\| r_{k}^{R} \|_{H}^{2}} - \frac{1}{\| r_{k-1}^{R} \|_{H}^{2}}} = \frac{\| r_{k}^{R} \|_{H}^{2}}{1- \frac{\| r_{k}^{R} \|_{H}^{2}}{\| r_{k-1}^{R} \|_{H}^{2}}},\\
  \| r_{k}^{W} \|_{H} &= \frac{\| r_{k}^{R} \|_{H}}{\sqrt{ 1- \frac{\| r_{k}^{R} \|_{H}^{2}}{\| r_{k-1}^{R} \|_{H}^{2}}}}
\cdot\end{align*}

We conclude that while Rapoport's method guarantees a minimal number of iterations through the monotonically non-increasing residual norm, the increased computational cost associated with obtaining the minimum residual may negate the potential savings in computation time.
Further, while the residual norms produced by Widlund's method are not monotonically non-increasing, we have shown that ``on average'' they are non-increasing.
Thus, we emphasize Widlund's method in our computational experiments.

    %! TEX root = **/000-main.tex
% vim: spell spelllang=en:

\section{Three Possible Implementations}%
\label{sec:currentsituation}

Widlund \cite{Widlund78} shows that in exact arithmetic, there are two possible formulations of Widlund's method. 
The first computes the basis vectors directly from a 3-term recurrence of the form
\begin{equation}
    \alpha_m v_{m+1} = -Kv_m - \beta_m v_{m-1} \label{eq:lanczos3termRecurr}
,\end{equation}
where $\alpha_{j+1} = 1+\frac{\rho_{j+1}}{\rho_{j}}\alpha_j$, $\rho_j = (v_j, v_j)_H$, $\beta_j = 1-\alpha_j$, and $\alpha_0 = 1$.
This gives rise to the Lanczos relation \eqref{eq:lanczosRelation}.
Then, with a special normalization, a solution can be formed as
\begin{equation}
    x_m = V_my_m, \text{ where } (I_m+T_m)y_m = e_1.
\end{equation}
Then the tridiagonal structure of $I_m + T_m$ is exploited to instead write 
\begin{equation}
    x_{m+1} = (1-\omega_{m+1})x_{m-1} + \omega_{m+1}x_m + \omega_{m+1}v_{m}, \quad \omega_{m+1} = 1/\alpha_m\label{eq:xRecurr}
.\end{equation}
This is rendered in Algorithm \ref{alg:cgwExpK}, which we refer to as ``ExplicitK''.
\begin{algorithm}[h]
         \begin{algorithmic}[1]  
             \vspace{0cm}
             \State $\alpha_0 = 0$, $x_0 = 0$, $v_{0} = 0$, $\rho_{0} = 1$
             \State $x_1 = 0$, $v_1 = H^{-1}b$, $r_{1} = b$
             \For{$m = 1,2,3,\dots$}
             \State $\rho_m = (v_m,v_m)_H$
             \State $\alpha_m = 1+\frac{\rho_m}{\rho_{m-1}}\alpha_{m-1}$, $\omega_{m+1} = 1/\alpha_m$
             \State $x_{m+1} = (1-\omega_{m+1})x_{m-1} + \omega_{m+1}x_m + \omega_{m+1}v_{m}$
             \State $r_{m+1} = b - Ax_{m+1}$
             \State $v_{m+1} = -\omega_{m+1}Kv_m + (1-\omega_{m+1}) v_{m-1}$
             \EndFor
         \end{algorithmic}
         \caption{Widlund: Implementation with Explicit $K$} \label{alg:cgwExpK}
 \end{algorithm}

 The second implementation, Algorithm \ref{alg:cgwHinv}, which was used in \cite{altmannThreetermRecurrenceIterations2026, Guducu.Liesen.Mehrmann.Szyld.22}, exploits the fact that \eqref{eq:xRecurr} and \eqref{eq:lanczos3termRecurr} result in
\begin{equation}
    v_{m} = H^{-1}r_m \label{eq:vEq}
,\end{equation}
and therefore 
\begin{equation}
    (v_{m}, v_{m})_H = (r_m, v_{m}) \label{eq:innerProdProperty}
.\end{equation}
We refer to Algorithm \ref{alg:cgwHinv} as ``implicitK''.
\begin{algorithm}[h]
         \begin{algorithmic}[1]  
             \vspace{0cm}
             \State $\alpha_0 = 0$, $x_0 = 0$, $\rho_{0}=1$
             \State $x_{1}=0$, $r_{1} = b$
             \For{$m = 1,2,3,\dots$}
             \State $v_{m} = H^{-1}r_m$
             \State $\rho_{m} = (v_{m}, r_m)$
             \State $\alpha_{m} = 1+\frac{\rho_{m}}{\rho_{m-1}}\alpha_{m-1}$, $\omega_{m+1} = 1/\alpha_{m}$
             \State $x_{m+1} = (1-\omega_{m+1})x_{m-1} + \omega_{m+1}x_m + \omega_{m+1}v_{m}$
             \State $r_{m+1} = b - Ax_{m+1}$
             \EndFor
         \end{algorithmic}
         \caption{Widlund: Implementation with implicit $K$ by solving $Hv_{m} = r_m$} \label{alg:cgwHinv}
 \end{algorithm}

 Moreover, one can show that \eqref{eq:vEq} and \eqref{eq:xRecurr} imply
 \begin{equation}
     r_{m+1} = (1-\omega_{m+1})r_{m-1}-\omega_{m+1}Sv_{m}, \label{eq:rRecurr}
 \end{equation}
resulting in a third possible algorithm, Algorithm \ref{alg:cgwHinvrRecurr}, which we refer to as ``3termR''.
\begin{algorithm}[h]
         \begin{algorithmic}[1]  
             \vspace{0cm}
             \State $\alpha_0 = 0$, $x_0 = 0$, $r_0 =b$, $\rho_{0} = 1$
             \State $x_{1} = 0$, $r_{1} = b$
             \For{$m = 1,2,3,\dots$}
             \State $v_{m} = H^{-1}r_m$
             \State $\rho_{m} = (v_{m}, r_m)$
             \State $\alpha_{m} = 1+\frac{\rho_{m}}{\rho_{m-1}}\alpha_{m-1}$, $\omega_{m+1} = 1/\alpha_{m}$
             \State $x_{m+1} = (1-\omega_{m+1})x_{m-1} + \omega_{m+1}x_m + \omega_{m+1}v_{m}$
             \State $r_{m+1} = (1-\omega_{m+1})r_{m-1}-\omega_{m+1}Sv_{m}$
             \EndFor
         \end{algorithmic}
         \caption{Widlund: Implementation by solving $Hv_{m} = r_m$, and $r_{m+1}$ recurrence} \label{alg:cgwHinvrRecurr}
 \end{algorithm}

% While the algorithms are mathematically equivalent in the case when $H^{-1}$ is applied exactly, this is no longer true when $H^{-1}$ is applied inexactly, which we show in section \ref{sec:inexAlgNotEquiv}.
We now explore how these three implementations behave in practice, together with two implementations available in the literature: the implementation in \cite{MehrmannStoll}, which we refer to as MSS, and \cite[Algorithm 2]{flexible-short-recurrence-Frommer}, which we refer to as DFK.

In MSS, \eqref{eq:lanczos3termRecurr} is used to construct an $(m+1)\times m$ tridiagonal matrix $T_{m+1,m}$, which is used in the auxiliary problem
\[
  \left( I + T_{m+1,m} \right)y_{m} = \| \hat{b} \|_{H} e_{1}
\cdot\] 
MSS also uses the following identity to simplify $T_{m+1,m}$.
\begin{align*}
  \alpha_{k} &= \left( \left( I + K \right)v_{m},v_{m-1} \right)_{H}\\
             &= - \left( v_{m},v_{m} + \alpha_{m-1}v_{m-1} + \beta_{m-2}v_{m-2} \right)_{H}\\
             &= -\beta_{m}
.\end{align*}
Finally, the new approximate solution is $x_{m} = V_{m} y_{m}$.
One can show that the residual, $r_m$, is a scalar multiple of $v_{m+1}$, where the scalar is determined by $y_{m}$ and $\beta_{m+1}$.

\begin{algorithm}
\caption{Widlund: MSS implementation}
\label{alg:MSS}
\begin{algorithmic}[1]
\State $\widehat b = H^{-1}b$, $\beta_0 = \|\widehat b\|_H$
\State $v_0 = 0$, $v_1 = \widehat b/\beta_0$, $\beta_1 = 0$
\State $\eta_1 = 1$, $\zeta_1 = \beta_0$
\For{$m = 1,2,3,\ldots$}
    \State $w = Kv_m + \beta_{m}v_{m-1}$
    \State $\widetilde{w} = Hw$
    \State $p_m = (v_m + \beta_m p_{m-1})/\eta_m$
    \State $x_m = x_{m-1} + \zeta_m p_m$
    \State $r_m = -(\zeta_m/\eta_m)\widetilde{w}$
    \State $\beta_{m+1} = (w,\widetilde{w})^{1/2}$
    \State $v_{m+1} = w/\beta_{m+1}$
    \State $\lambda_{m+1} = \beta_{m+1}/\eta_m$
    \State $\eta_{m+1} = 1 + \lambda_{m+1}\beta_{m+1}$
    \State $\zeta_{m+1} = -\lambda_{m+1}\zeta_m$
\EndFor
\end{algorithmic}
\end{algorithm}

DFK uses right preconditioning to form the Lanczos relation
\[
\left( H + S \right)Z_{m} = V_{m+1}T_{m+1,m}
,\] 
where the columns of $V_{m+1}$ are $H^{-1}$-orthonormal and form a basis for the Krylov subspace $\mathcal{K}\left( I + SH^{-1},b \right)$, and $Z_m = H^{-1}V_m$.
The residuals can be written analogous to the residual computation used in MSS.
That is, $r_m$ is expressed as a scalar multiple of $v_{m+1}$, where the scalar is determined by $y_{m}$ and the entries of $T_{m+1,m}$.

\begin{algorithm}[h]
\label{alg:DFK}
\begin{algorithmic}[1]
\State $x_0=0$, $r_0=b$
\State $\widehat{r}_0 = H^{-1}r_0$, $\beta_0=(r_0,\widehat{r}_0)^{1/2}$
\State $v_1=r_0/\beta_0$, $z_1=\widehat{r}_0/\beta_0$, $v_0=0$, $z_0=0$
\For{$m=1,2,3,\ldots$}
    \State $w=(H+S)z_m$
    \State $\alpha_m=(w,z_m)$
    \State $\gamma_m=(w,z_{m-1})$
    \State $w=w-\alpha_m v_m-\gamma_m v_{m-1}$
    \State $\widehat{w}= H^{-1}w$
    \State $\beta_m=(w,\widehat{w})^{1/2}$
    \State $v_{m+1}=w/\beta_m$
    \State $z_{m+1}=\widehat{w}/\beta_m$
    \State Solve $T_m y_m = \beta_0 e_1$
    \State $x_m = Z_m y_m$
    \State $\eta_m = e_m^T y_m$
    \State $r_m = -\beta_m\eta_m v_{m+1}$
\EndFor
\end{algorithmic}
\caption{Widlund: DFK implementation}
\end{algorithm}

The numerical examples in this section and the next, following \cite{convection-diffusion-example}, are based on the advection-diffusion equation for a concentration $c = c\left( x,y,t \right): \left[ 0,1 \right] \times \left[ 0,1 \right] \times \left[ 0,T \right] \to \mathbb{C}$, where $0 < T \in \mathbb{R}$,
\begin{align}
  \frac{\partial}{\partial t} c = &-\frac{\partial}{\partial x} \alpha \frac{\partial}{\partial x}c - \frac{\partial}{\partial y}\beta \frac{\partial}{\partial y}c + \frac{\partial}{\partial x}\nu c + \frac{\partial}{\partial y}\mu c,\\
                                &c\left( x,y,t \right) = 0 \text{ if } \left( x,y \right) \in \partial \left( \left[ 0,1 \right] \times \left[ 0,1 \right] \right),\\
                                &c(x,y,0) = c_{0} \left( x,y \right) \text{ if } \left( x,y \right) \in \left[ 0,1 \right] \times \left[ 0,1 \right] 
.\end{align}
Both the positive diffusion coefficients $\alpha,\beta$ and the velocity field $\left( \nu,\mu \right)$ may depend on $x$ and $y$.
The discretization in space uses the standard central finite difference scheme with $N$ equidistant points in both $x$ and $y$.
This produces a semi-discretized system
\begin{equation}
  \frac{\partial}{\partial t} c = B t, \quad c\left( x,y,0 \right) = c _{0}\left( x,y \right) \text{ for all grid points }x,y
,\end{equation}
where $c = c\left( t \right)$ is a two-dimensional array with elements corresponding to the grid points.
The implicit Euler rule is used as a symplectic integrator, so for a given time $t$, step size $k$, and $A = I + \frac{k}{2} B$, we have to solve
 \begin{equation}
   A c\left( t+1 \right) = \left( I + \frac{k}{2} B \right)c \left( t+1 \right) = \frac{k}{2}B c\left( t \right)\label{eq:numericalExampleEuler}
\end{equation}
for $c\left( t+1 \right)$.
We set the exact solution, $x$, to be a random vector whose components are independent and uniformly distributed on $\left( 0,1 \right)$ with corresponding right-hand side $b = Ax$.

For the first example, we take $N = 100$ in $x$ and $y$, $h = \frac{1}{N+1}$, and $k = 0.5h^{2}$.
Thus, $n=10,000$.
We consider the case of advection in one direction, $\nu = 50$, $\mu = 0$, and constant diffusion, $\alpha = \beta = 0.5$.

For the second example, following the example in \cite{Widlund78}, we take $N = 127$, $\alpha = \beta = 1$, $\nu = 10000$, and $\frac{\partial}{\partial t}c = 0$.
Thus, $n=16,129$.
We take the exact solution and right-hand side to be the same as in example one.

\begin{figure}[htb]
    \centering
    \includegraphics[width=\textwidth]{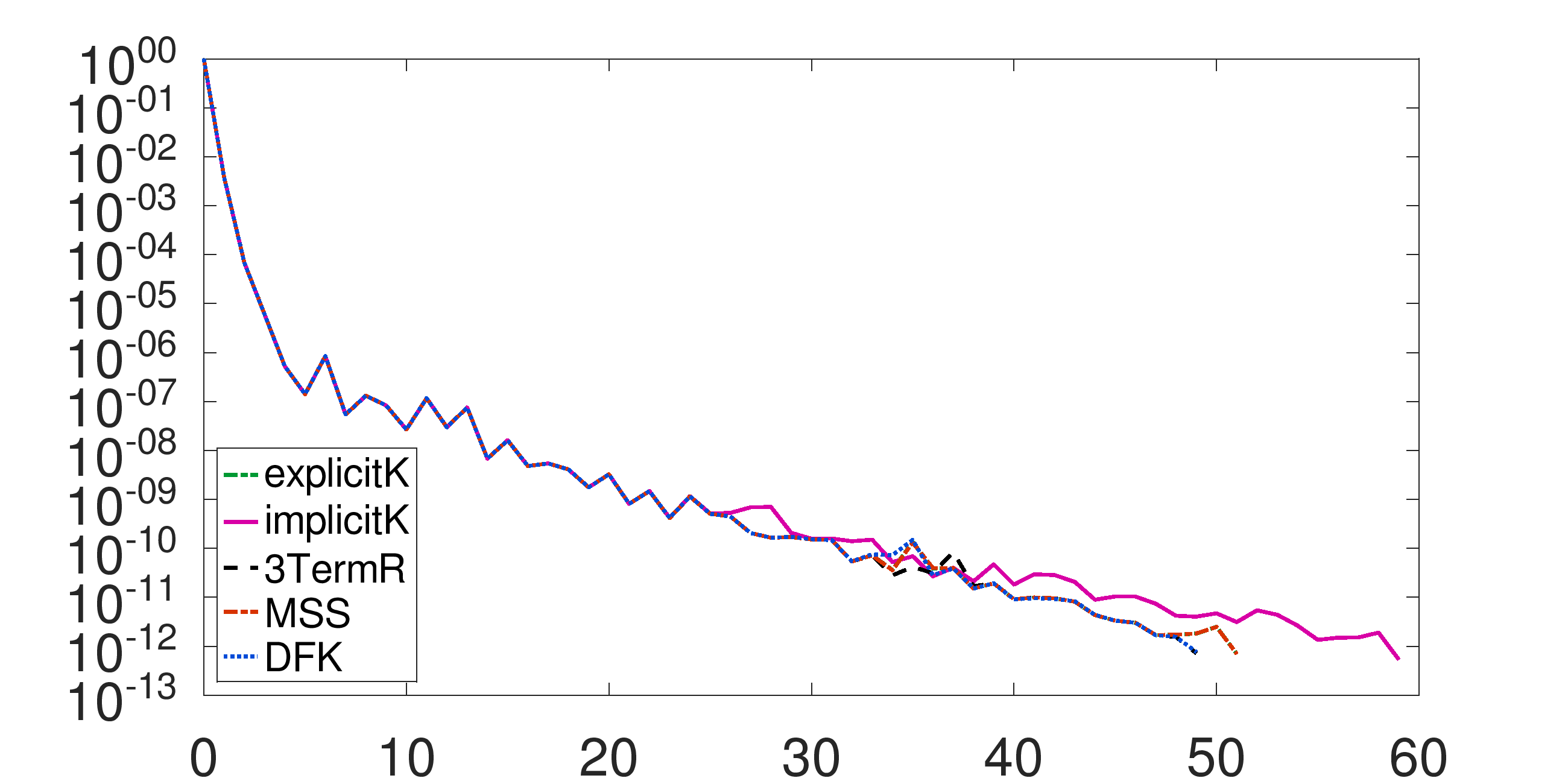}
    \caption{$H$-norm relative residual for example one ($n=10000$).}
    \label{fig:prob1ResExact}
\end{figure}

In Figure \ref{fig:prob1ResExact}, we report on experiments with the five different implementations for example one.
We observe that the residual curves for all five implementations largely coincide, with 3termR converging in the fewest iterations and implicitK requiring the most iterations.

\begin{figure}[htb]
    \centering
    \includegraphics[width=\textwidth]{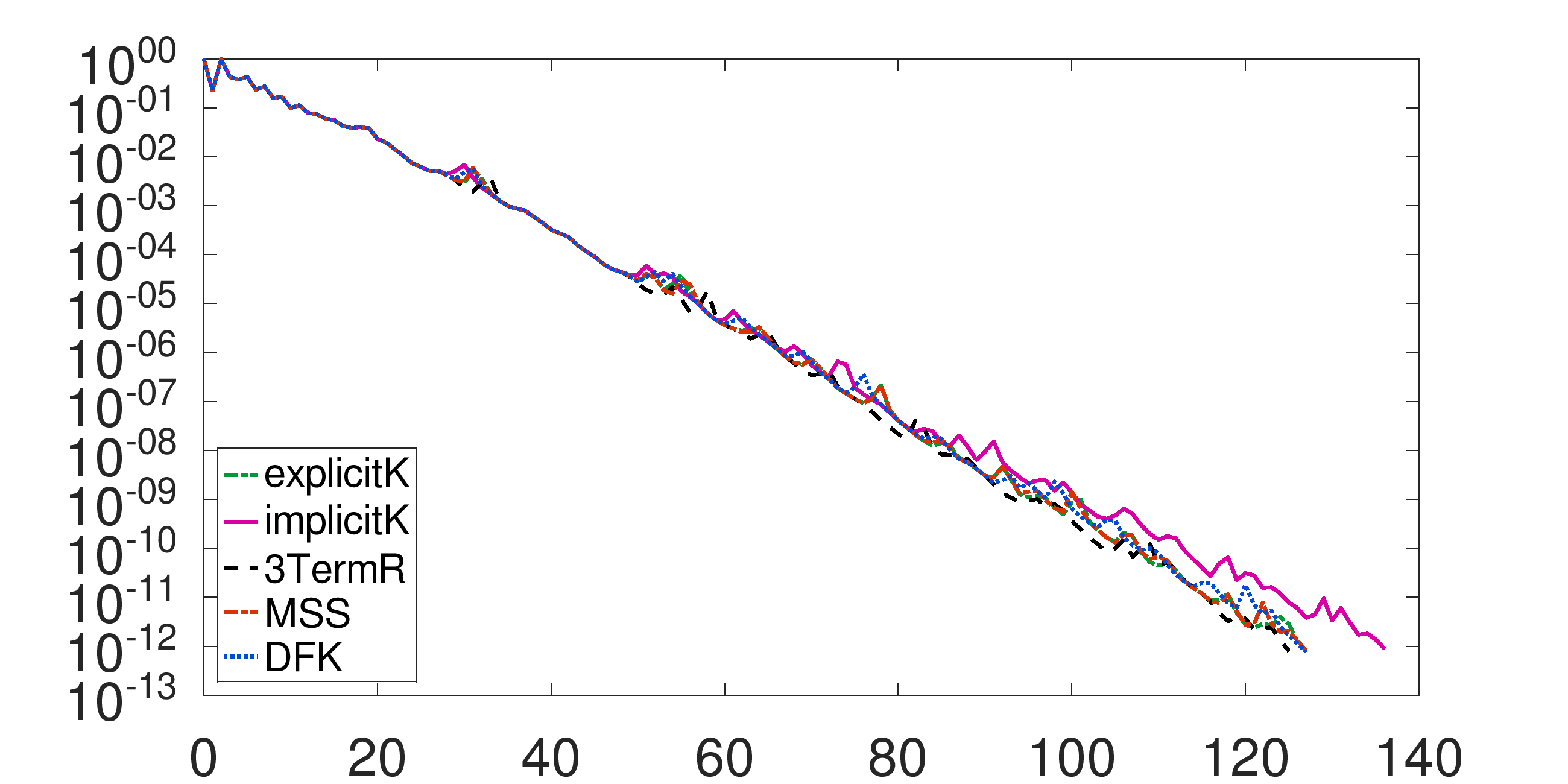}
    \caption{$H$-norm relative residual for example two ($n=16129$).}
    \label{fig:prob2ResExact}
\end{figure}

In Figure \ref{fig:prob2ResExact}, we report on experiments with the five different implementations for example two.
We observe that the residual curves remain similar overall, with 3termR again converging fastest and implicitK slowest.

Since the number of iterations required for each implementation is similar, we will use the cost per iteration, both in storage and computation, as the main metric for recommending the implementation.
All 5 implementations require one solve.
ExplicitK, ImplicitK, MSS, and 3termR all store 6 vectors, whereas DFK stores 10 vectors.
ExplicitK requires three matvecs, namely $Av$, $Hv$, and $Sv$.
MSS requires two matvecs, namely $Hv$ and $Sv$. 
ImplicitK and DFK both require one matvec, $Av$.
Finally, 3termR also only requires one matvec, $Sv$.
While 3termR does not compute the true residual each iteration, the max distance from the computed residual to the true residual in both examples is on the order of $10^{-15}$.
The following table contains the number of iterations, as well as the time per iteration, for each implementation applied to example two.

\begin{table}[ht]
\centering
\begin{tabular}{lcccc}
\hline
 & Iters & ms/iter & Total & vs.\ Best \\
\hline
\textbf{3termR}
 & \textbf{125} & \textbf{3.253} & \textbf{406.6 ms} & --- \\
ImplicitK
    & 137 & 3.600 & 493.3 ms & $1.21\times$ \\
DFK
    & 127 & 4.179 & 530.7 ms & $1.31\times$ \\
MSS
    & 127 & 4.269 & 542.2 ms & $1.33\times$ \\
ExplicitK
    & 127 & 4.462 & 566.6 ms & $1.39\times$ \\
Rapoport
    & 126 & 4.549 & 573.2 ms & $1.41\times$ \\
\hline
\end{tabular}
\end{table}

Overall, the computational study indicates that 3termR requires the fewest iterations, the fewest matvecs, and uses the least amount of storage.
We recommend 3termR in the exact setting.

    %! TEX root = **/000-main.tex
% vim: spell spelllang=en:

\section{Inexact Widlund Implementations}%
\label{sec:inexAlgs}

In practical situations, we may not solve for $r_{m}$ exactly in step 4 of implicitK or 3termR.
For example, we may have an inexact $LU$ factorization of $H$, and we solve
\[
LUv_{m} = r_{m}
\] 
as done, e.g., in \cite{Guducu.Liesen.Mehrmann.Szyld.22}.
Thus, the basis vectors are no longer $H$-orthogonal.
Nevertheless, the implementations can work.
There are two possible situations.
First, we have a fixed $\bar{H}$ approximation to $H$ for all iterations, as in the above example of inexact $LU$. 
Second, we have $\bar{H} = \bar{H}_{m}$ changing with each iteration, e.g., when the system $Hv_{m} = r_{m}$ is solved with an iterative method to a specified tolerance.
We explore the first case computationally in the remainder of this section.

To define the inexact analogues of the implicitK and 3termR algorithms, which we refer to as ``inexactImplicitK'' and ``inexact3termR,'' we replace line 4 of both algorithms by
\[
v_{m} = \bar{H}^{-1}r_{m}
\cdot\] 
We refer to \cite[Section 3.2]{flexible-short-recurrence-Frommer} for the inexact version of DFK, which we refer to as ``implicitDFK''.
Our examples are the same as in Section \ref{sec:currentsituation}, approximating the action of $H^{-1}$ using either incomplete Cholesky factorization with threshold $10^{-4}$ or incomplete Cholesky with no fill.
Unless otherwise specified, we use the threshold incomplete Cholesky.

\begin{figure}[htb]
    \centering
    \includegraphics[width=\textwidth]{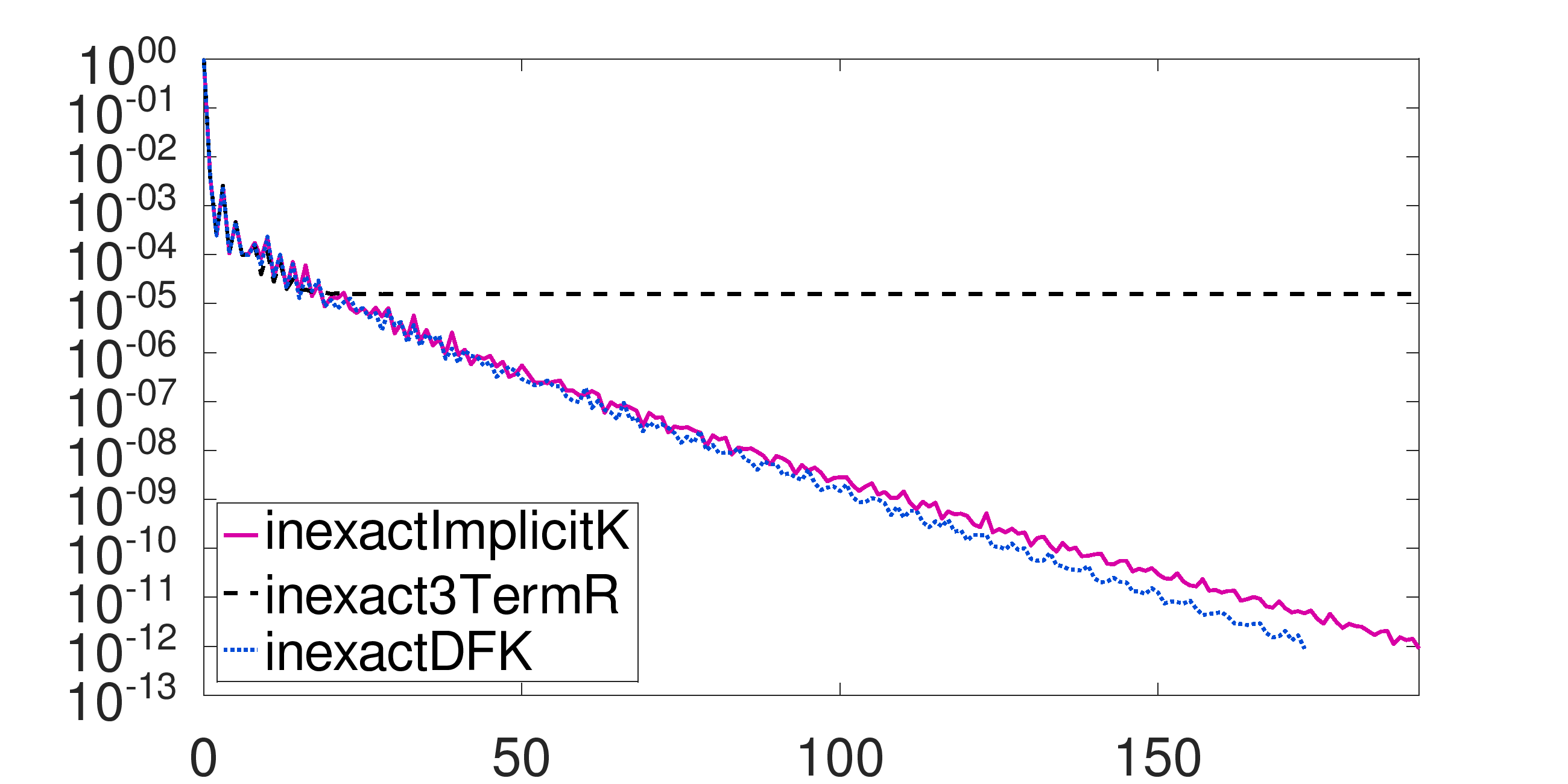}
    \caption{$H$-norm relative residual for example one}
    \label{fig:prob1Res}
\end{figure}

In Figure \ref{fig:prob1Res}, we report on experiments with inexactImplicitK, inexact3termR, and inexactDFK for example one.
We observe that inexact3termR stagnates, while inexactImplicitK and inexactDFK converge and have similar residual curves.
Because inexact3termR stagnates, we do not recommend it.
Also, it is worth noting that the number of iterations needed for inexactImplicitK and inexactDFK to converge is more than three times as large as in the exact case.
It is notable that inexactImplicitK converges at all, since the inexactness causes the basis generated from implicitK to lose the orthogonality needed in the three-term recurrence of the basis vectors.
However, since the level of inexactness is relatively small given our method of doing incomplete Cholesky, \cite[Proposition 5]{flexible-short-recurrence-Frommer} provides a guarantee of convergence for inexactDFK.

\begin{figure}[htb]
    \centering
    \includegraphics[width=\textwidth]{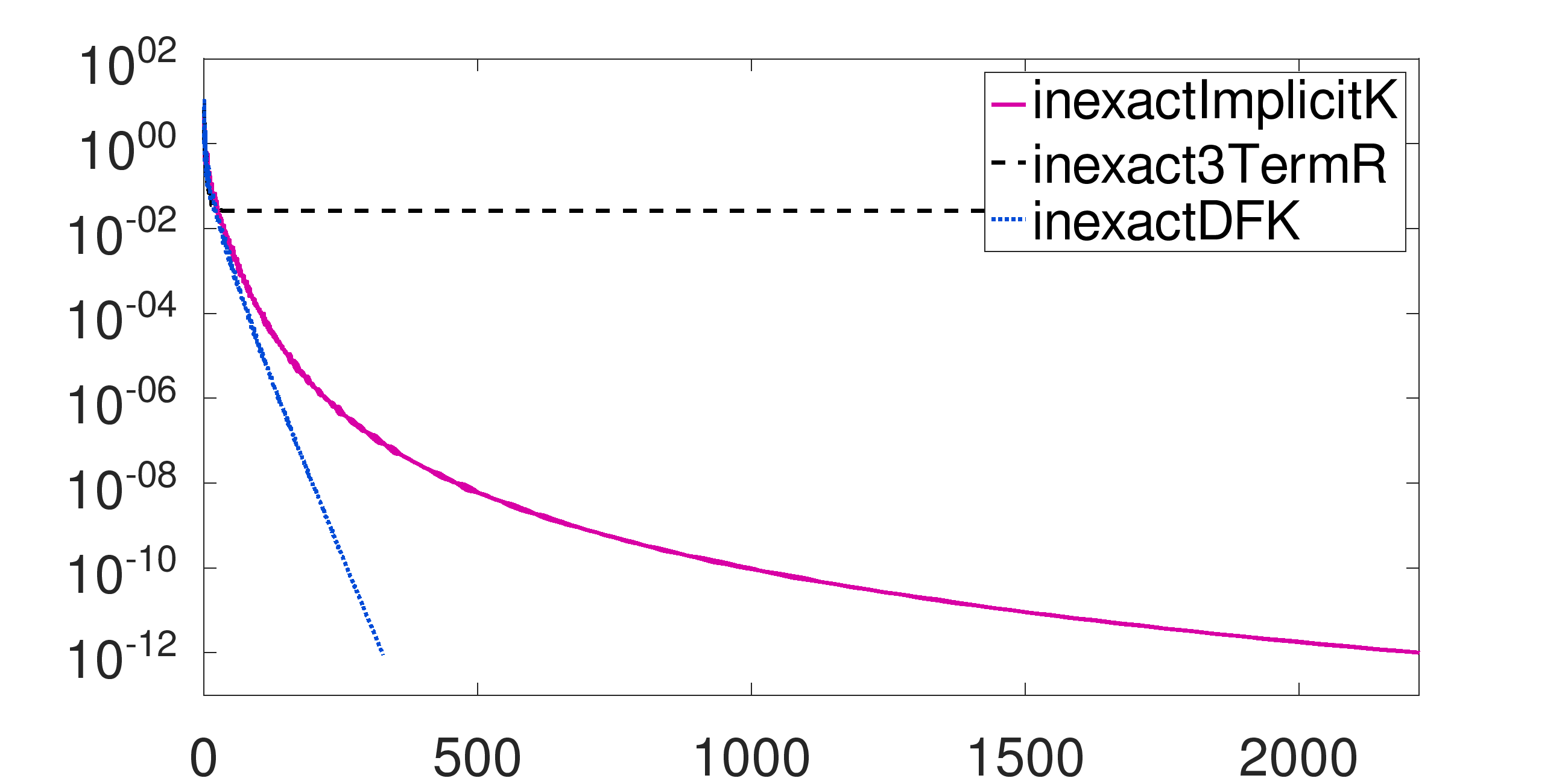}
    \caption{ $H$-norm relative residual for example two.}
    \label{fig:prob2ResNotRandomZ}
\end{figure}

% \begin{figure}[htb]
%     \centering
%     \includegraphics[width=\textwidth]{prob2Res.png}
%     \caption{ $H$-norm relative residual for example two.}
%     \label{fig:prob2Res}
% \end{figure}

In Figure \ref{fig:prob2ResNotRandomZ}, we report on experiments with the three inexact implementations for example two.
We observe that inexact3termR stagnates again, and that both inexactDFK and inexactImplicitK converge.
However, in this example inexactImplicitK takes significantly longer to converge than inexactDFK.

\begin{figure}[htb]
    \centering
    \includegraphics[width=\textwidth]{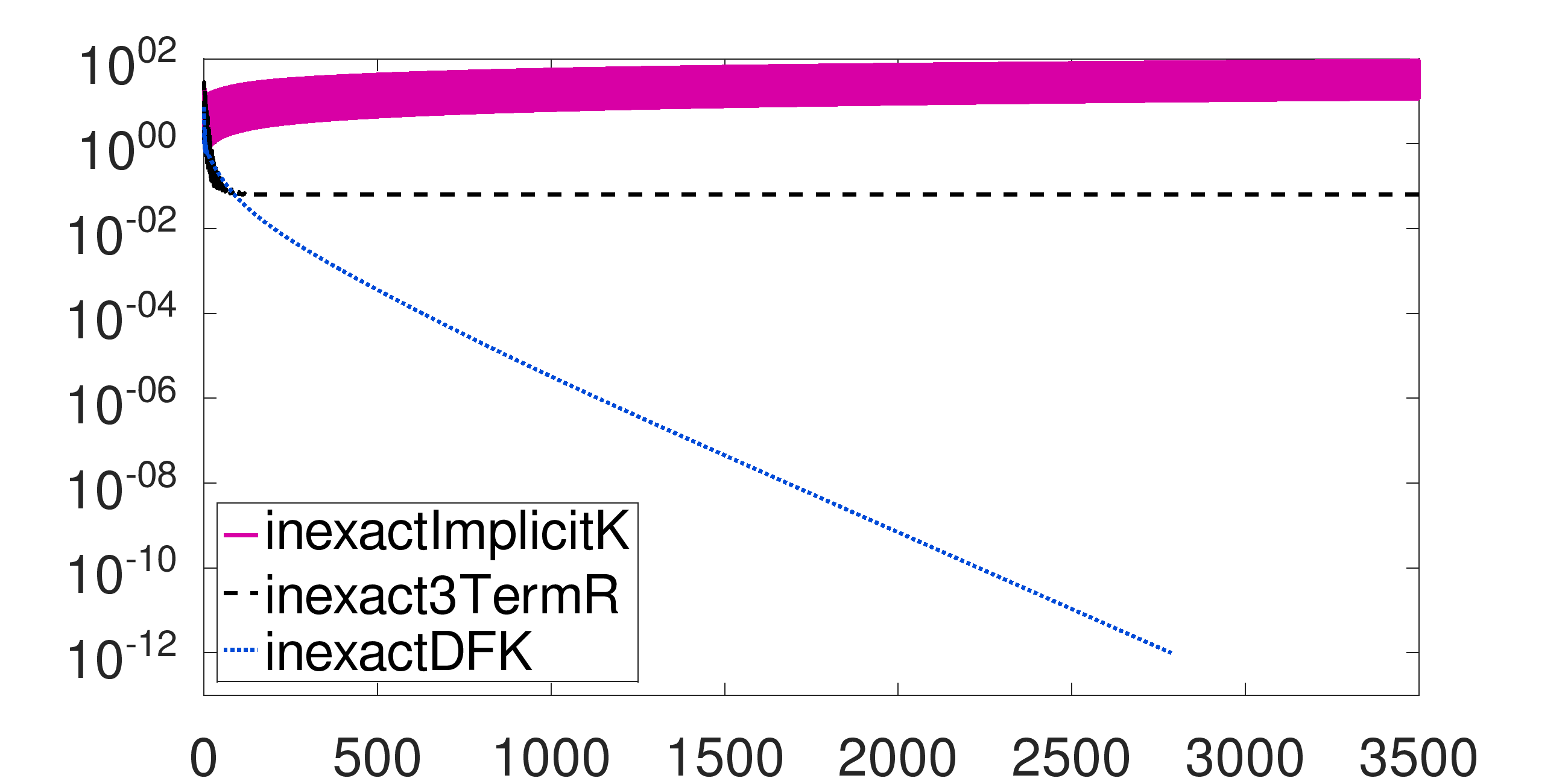}
    \caption{ $H$-norm relative residual for example two with no fill. }
    \label{fig:prob2ResNoFill}
\end{figure}

In Figure \ref{fig:prob2ResNoFill} we report on experiments with the three inexact implementations for example two when we use incomplete Cholesky with no fill.
While inexact3termR and inexactDFK behave similarly to the experiment shown in Figure \ref{fig:prob2ResNotRandomZ}, inexactImplicitK actually diverges.
More specifically, the even and odd residual norms of inexactImplicitK both stagnate, though to different values.

\begin{figure}[htb]
    \centering
    \includegraphics[width=\textwidth]{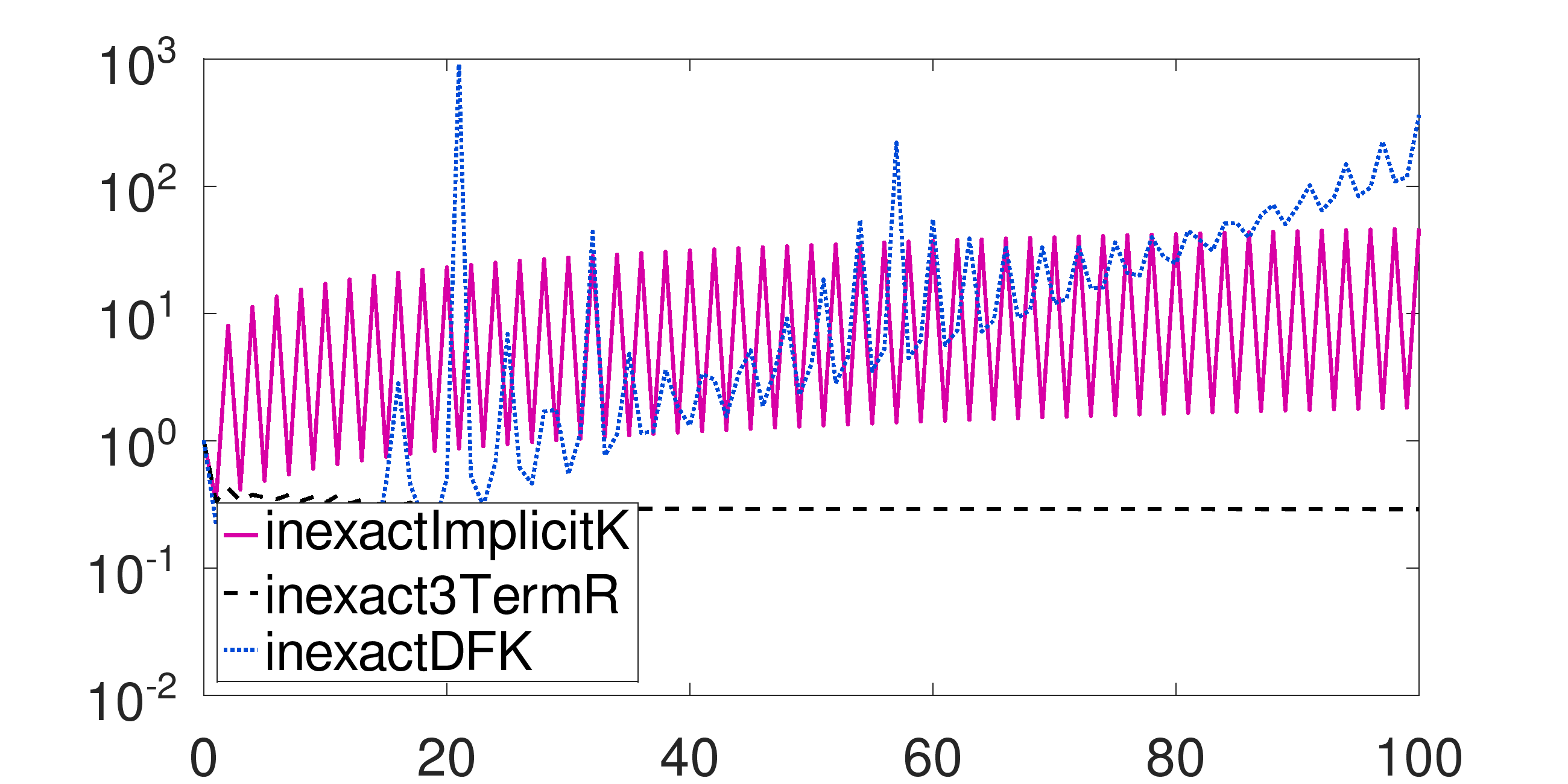}
    \caption{ $H$-norm relative residual for Stokes problem with no fill. }
    \label{fig:probStokesResNoFill}
\end{figure}

In Figure \ref{fig:probStokesResNoFill} we report on experiments with the three inexact implementations for the Stokes problem setup outlined in Section \ref{sec:nearOptimal}, and we compute incomplete Cholesky with no fill.
We observe that inexact3termR stagnates as usual, and inexactImplicitK also displays the same long term behaviour as in the experiment in Figure \ref{fig:prob2ResNoFill}.
However, inexactDFK now diverges due to the level of inexactness of the factorization.
Table~\ref{tbl:example1} and Table \ref{tbl:example2} show the number of iterations and wall-clock time needed for both methods to converge for a range of drop tolerances used in the incomplete Cholesky factorization of $H$.
The lines in the ``inexactImplicitK'' column indicate that the implementation did not converge within 3500 iterations.

\begin{table}[htbp]
\centering
\caption{\textbf{Example 1} ($n = 10 000$).}
\begin{tabular}{r r r r r r}
\hline
 & & \multicolumn{2}{c}{inexactImplicitK} & \multicolumn{2}{c}{inexactDFK} \\
droptol & nnz(L) & Iters & Time (s) & Iters & Time (s) \\
\hline
no-fill & 29 800 & \multicolumn{2}{c}{---} & \textbf{2873} & 1.57 \\
1e-01 & 29 800 & \multicolumn{2}{c}{---} & \textbf{102} & 0.05 \\
1e-02 & 49 303 & \multicolumn{2}{c}{---} & \textbf{84} & 0.05 \\
1e-03 & 123 438 & \multicolumn{2}{c}{---} & \textbf{162} & 0.11 \\
1e-04 & 276 914 & \textbf{785} & 0.68 & \textbf{277} & 0.27 \\
1e-05 & 536 252 & \textbf{305} & 0.44 & \textbf{270} & 0.43 \\
1e-06 & 792 100 & \textbf{203} & 0.41 & \textbf{198} & 0.42 \\
1e-07 & 898 097 & \textbf{172} & 0.41 & \textbf{168} & 0.42 \\
1e-08 & 931 242 & \textbf{156} & 0.39 & \textbf{154} & 0.40 \\
1e-09 & 947 565 & \textbf{149} & 0.38 & \textbf{147} & 0.40 \\
1e-10 & 957 773 & \textbf{145} & 0.39 & \textbf{141} & 0.39 \\
exact & 206 332 & \textbf{135} & 0.10 & \textbf{125} & 0.11 \\
\hline
\end{tabular}
\label{tbl:example1}
\end{table}

\begin{table}[htbp]
\centering
\caption{\textbf{Example 2} ($n = 16 129$).}
\begin{tabular}{r r r r r r}
\hline
 & & \multicolumn{2}{c}{inexactImplicitK} & \multicolumn{2}{c}{inexactDFK} \\
droptol & nnz(L) & Iters & Time (s) & Iters & Time (s) \\
\hline
no-fill & 48 133 & \multicolumn{2}{c}{---} & \textbf{3222} & 4.45 \\
1e-01 & 48 133 & \multicolumn{2}{c}{---} & \textbf{112} & 0.16 \\
1e-02 & 79 759 & \multicolumn{2}{c}{---} & \textbf{128} & 0.19 \\
1e-03 & 201 306 & \multicolumn{2}{c}{---} & \textbf{216} & 0.42 \\
1e-04 & 460 730 & \multicolumn{2}{c}{---} & \textbf{353} & 1.11 \\
1e-05 & 935 444 & \multicolumn{2}{c}{---} & \textbf{351} & 1.86 \\
1e-06 & 1 504 060 & \textbf{242} & 2.04 & \textbf{228} & 2.00 \\
1e-07 & 1 798 165 & \textbf{187} & 1.87 & \textbf{180} & 1.91 \\
1e-08 & 1 888 277 & \textbf{165} & 1.70 & \textbf{165} & 1.78 \\
1e-09 & 1 927 546 & \textbf{155} & 1.54 & \textbf{154} & 1.59 \\
1e-10 & 1 951 080 & \textbf{149} & 1.60 & \textbf{147} & 1.65 \\
exact & 350 112 & \textbf{137} & 0.37 & \textbf{127} & 0.39 \\
\hline
\end{tabular}
\label{tbl:example2}
\end{table}

From our experiments, we recommend inexactDFK in the inexact setting. 
While it requires four more vectors of storage and the computational cost per iteration is higher than inexactImplicitK, the convergence of inexactDFK is faster when both methods converge which causes the total cost to be similar in most cases.
InexactDFK is also more robust to the level of inexactness and has more theoretical guarantees,~see~\cite{flexible-short-recurrence-Frommer}.

We remark that inexact3termR can be interpreted as the exact version, 3termR, applied to $\bar{A}x = (\bar{H} + S)x = b$.
This can be shown inductively using 
\begin{align*}
b-\bar{A}x_{m+1}
&= b-\bar{A}
\left[
  (1-\omega_{m+1})x_{m-1}
  +\omega_{m+1} x_m
  +\omega_{m+1} v_m
\right] \\
&= (1-\omega_{m+1})\left(b-\bar{A}x_{m-1}\right)
+\omega_{m+1}\left(b-\bar{A}x_m\right)
-\omega_{m+1}\bar{A}v_m \\
&= (1-\omega_{m+1})r_{m-1}
+\omega_{m+1} r_m
-\omega_{m+1}\bar{H}v_m
-\omega_{m+1} Sv_m \\
&= (1-\omega_{m+1})r_{m-1}
-\omega_{m+1} Sv_m \\
&= r_{m+1}.
\end{align*}

The iterates of inexact3termR converge to $\bar{x}_*$, the true solution of $\bar{A}x = b$.
Thus, the residual norm converges to $\|(I - A \bar{A}^{-1})b\|_H$, which is not necessarily close to zero.

    % \input{115-inexact-not-eq}
    % \input{116-inexact-props}
    %! TEX root = **/000-main.tex
% vim: spell spelllang=en:

\section{Conclusion}%
\label{sec:conclusion}

We proved that the residual norms of Widlund's method, while they are not monotonically decreasing, are ``nearly-optimal'' in a certain sense; and this result justifies its practical usage over Rapoport's method.
Through a computational study, we demonstrated that 3termR should be used as the implementation in the exact setting, and inexactDFK should be used in the inexact setting.

%%%%%%%%%%%%%%%%%%%%%%%%%%%%%%%%%%%%%%%%%%%%%%%%%%%%%%%%%%%%%%%%%%%%%%%%%%%%%%%%
% BIBLIO
%%%%%%%%%%%%%%%%%%%%%%%%%%%%%%%%%%%%%%%%%%%%%%%%%%%%%%%%%%%%%%%%%%%%%%%%%%%%%%%%

    % \printbibliography

    \bibliography{biblio.bib}

%%%%%%%%%%%%%%%%%%%%%%%%%%%%%%%%%%%%%%%%%%%%%%%%%%%%%%%%%%%%%%%%%%%%%%%%%%%%%%%%
% APPENDIX
%%%%%%%%%%%%%%%%%%%%%%%%%%%%%%%%%%%%%%%%%%%%%%%%%%%%%%%%%%%%%%%%%%%%%%%%%%%%%%%%

    % uncomment to add appendix section prefix to numbering
    %\appendixpagenumbering

    % \appendix

    % \include{990-appendix}

\end{document}